\documentclass[a4paper, 10pt]{paper}
\usepackage{geometry}
\usepackage[english]{babel}

\usepackage[utf8x]{inputenc}
\usepackage{amsmath, amsthm, amssymb}
\usepackage{mathtools}
\usepackage{enumitem}
\usepackage{nameref,hyperref,cleveref}
\usepackage{thmtools}
\usepackage{dsfont,bbm}
\usepackage{graphicx}
\usepackage{xcolor}
\usepackage{textcomp}
\usepackage{comment}

\usepackage{wrapfig}
\usepackage{pdfpages}
\usepackage{comment}

\title{Existence of a plane without edge crossings in projections of the random geometric graph}
\author{Lianne de Jonge\thanks{l.de.jonge.math@proton.me}\and Kinga Nagy \thanks{Osnabrück University, Germany, kinga.nagy@uni-osnabrueck.de} \thanks{Corresponding author}}
\date{\today}

\theoremstyle{definition}

\newtheorem*{remark*}{Remark}

\theoremstyle{plain}
\newtheorem{theorem}{Theorem}
\newtheorem{lemma}{Lemma}
\newtheorem*{corollary*}{Corollary}

\newtheorem{claim}{Claim}
\newtheorem{proposition}{Proposition}
\newtheorem*{claim*}{Claim}

\newcommand{\EE}{\mathbb{E}}
\newcommand{\PP}{\mathbb{P}}
\newcommand{\VV}{\mathbb{V}}
\newcommand{\Poi}{\mathrm{Poi}}
\newcommand{\ind}{\mathbbm 1}
\newcommand{\Cov}{\mathrm{Cov}}

\DeclareMathOperator{\dTV}{d_{TV}}

\DeclareMathOperator{\diam}{diam}

\renewcommand{\l}{\left}
\renewcommand{\r}{\right}

\newcommand{\dd}{{\rm d}}

\newcommand{\intl}{\int\limits}

\newcommand{\cross}{\mathrm{cr}}
\DeclareMathOperator{\spanvec}{span}

\newcommand{\radt}{r_t}

\newcommand{\R}{\mathbb{R}}
\newcommand{\N}{\mathbb{N}}

\renewcommand{\c}{\mathcal}
\renewcommand{\bf}{\mathbf}

\begin{document}

\maketitle
\begin{abstract}
    Consider a random geometric graph $G$ with a vertex set defined by a Poisson point process  with intensity $t>0$ in a convex body. We can generate a drawing of the graph by projecting the construction onto some plane $L$. Choosing different planes leads to different drawings, and in particular, potentially more or fewer edge crossings. In this paper, we prove that if the connection radius is smaller than a given threshold, the probability that there exists a plane with zero crossings tends to one as $t\to \infty$. We also state the asymptotic probability that such a plane is found after considering a given number of randomly chosen planes.
\end{abstract}

\noindent \textsc{Keywords.} graph crossing number, random geometric graph, Poisson point process, limit theorems. \\
\noindent \textsc{MSC classification.} 
60F05, 
60D05, 
68R10. 

\section{Introduction}\label{sec:intro}
To effectively communicate the structure of a graph $G=(V,E)$, we often create a drawing in which the vertices are represented by points and the edges by line segments. Depending on where we place the vertices, the line segments might need to cross. Since graph drawings with few crossings tend to be the most aesthetically pleasing, we want to minimize the number of crossings. The minimal number of crossings is called the graph crossing number and finding this number is a classical problem in graph theory, see for example the survey \cite{schaefer2013}.

Finding the crossing number is a very computationally expensive problem. Instead of attempting to create a perfect drawing, we could therefore also try to find an algorithm to create a drawing with few crossings: in this paper, we consider graph drawings of a random geometric graph (RGG), also referred to as the Gilbert graph, created by projecting the graph onto a plane.

The RGG is generated using a Poisson points process of intensity $t$ in a convex body $W\subset \R^d$, which defines the vertices of the graph. Two distinct vertices $v$ and $w$ are connected by an edge whenever $\|v-w\|\leq \radt$, where $\|\cdot\|$ denotes the Euclidean norm and $r_t>0$ is a given distance parameter. Using orthogonal projection onto a plane, we create a drawing of the resulting graph. This algorithm for graph drawing was introduced in \cite{Chimani2018}, where the expectation and the variance of the number of crossings in the drawing were derived. In particular, it was shown that the expected crossing number is of the order $t^4 r_t^{2d+2}$. 
The paper \cite{Döring2025} builds on these results and shows that, if $t^4 r_t^{2d+2}\to c\in(0,\infty)$, the point process created by the crossings on the plane converges to a Poisson point process. 

To find a nice drawing, we could try different projection planes. In the best-case scenario, we would find a plane such that no edges cross in the projection. In this paper, we show that there exists a constant $c^*$ such that if $t^4 r_t^{2d+2} \leq  c'\log t$ for some $c'<c^*$, the probability that such a plane exists converges to one as $t$ tends to infinity.

The existence of a plane without crossings does not guarantee that we can realistically find one. In our second result, we quantify the likelihood of finding such a plane, and in particular show that if $t^4 r_t^{2d+2}\to c \in(0,\infty)$, this probability is close to one after enough attempts.

This paper is an extension of the earlier arXiv version 
written by the first author, where only the three-dimensional case was considered, and using slightly different tools. That preprint also appeared as part of the first author's PhD thesis \cite{deJongeThesis}.

\subsection{Definitions and results}
Let $W\subset \R^d$ be a convex body (that is, a compact convex set with nonempty interior) for some $d\geq 3$, and $\eta_t$ a homogeneous Poisson point process on $W$ with intensity $t>0$.
The RGG is the graph $G=(V,E)$ with vertex set $V=\eta_t$ and edge set $E=\{\{v,w\}\subset \eta_t\colon\, 0<\|v-w\|\leq r_t\}$ for some connectivity radius $r_t>0$ that converges to zero as $t$ goes to infinity. The line segment between two points $v$ and $w$ is denoted by $[v,w]$. 
We consider orthogonal projections of the graph onto a plane $L$; in particular, we may assume that $L$ is a two-dimensional linear subspace. The space of two-dimensional linear subspaces of $\R^d$, called the Grassmanian and denoted by $\bf G_2(d)$, is equipped with a uniform Haar measure, denoted by $\nu_d$. Note that $\nu_d$ is a probability measure, that is, $\nu_d(\bf G_2(d))=1$.
We denote the projection of a point $v\in\R^d$ and a set $A\subset \R^d$ onto $L$ by $v|_L$ and $A|_L$ respectively. 

For a plane $L\in \bf G_2(d)$, define the function $h_L:(\R^d)^4\to \R$,
\begin{align}
\begin{split}\label{eq:kernel}
    h_L(v_1,v_2,v_3,v_4) 
    \coloneq & \ind \l(\|v_1-v_2\|\leq \radt,\, \|v_3-v_4\|\leq \radt,\, [v_1,v_2]|_{L}\cap [v_3,v_4]|_{L}\neq \emptyset\r)\\
    + & \ind \l( \|v_1-v_3\|\leq \radt,\, \|v_2-v_4\|\leq \radt,\, [v_1,v_3]|_{L}\cap [v_2,v_4]|_{L}\neq \emptyset\r)\\
    + & \ind \l(\|v_1-v_4\|\leq \radt,\, \|v_2-v_3\|\leq \radt,\, [v_1,v_4]|_{L}\cap [v_2,v_3]|_{L}\neq \emptyset\r).
\end{split}
\end{align}
For distinct points, this is the indicator that a $4$-tuple of vertices forms a crossing in the projection onto $L$; note that any tuple of distinct points can form at most one crossing, i.e. at most one of the indicators is $1$, and thus $h_L$ is indeed an indicator. 
The number of crossings on the plane $L$, denoted by $\cross(L)$, is then given by the Poisson $U$-statistic with symmetric kernel function $h_L/24$:
\[\cross(L)\coloneq \frac{1}{24}\sum_{(v_1,v_2,v_3,v_4)\in \eta_{t,\neq}^4} h_L(v_1,v_2,v_3,v_4),\]
where $\eta_{t,\neq}^4$ denotes the distinct $4$-tuples of $\eta_t$, and the prefactor $1/24$ appears due to the multiple counts.

Our first result is about the existence of a plane without crossings.

\begin{theorem}[Existence of a plane without crossings]\label{thm:existence}
    Let $X_t$ be the measure of planes that contain zero crossings:
    \[X_t\coloneqq \nu_d\big(\l\{L\in \bf G_2(d)\colon \cross(L)=0\r\}\big).\]
    Then, there exists $c^*=c^*(d,W)>0$ 
    such that if $c'<c^*$ and
    \[t^4 r_t^{2d+2}\leq c' \ln t\]
    for large enough $t$, then
    $$\lim_{t\to\infty}\PP\l(X_t>0\r) = 1.$$
    \label{thm:existence plane}
\end{theorem}

Note that from \cite[Theorem~1]{Chimani2018}, the expected number of crossings on a fixed plane $L$ is of order $t^4 r_t^{2d+2}$. As a consequence, \Cref{thm:existence} states that if the expected number of crossings on a fixed plane is at most some specified multiple of $\ln t$, there must still exist some plane without crossings with high probability.

In general, the almost sure existence of a plane with zero crossings does not imply that we can find such a plane within reasonably many tries. Our next result deals with this question.

\begin{theorem}[Finding a plane without crossings]\label{thm:finding}
    Let $L_1,L_2,\ldots$ be independent, uniformly distributed planes in $\bf G_2(d)$, and $m\in \N$ arbitrary. 
    Then the following hold:
    \begin{enumerate}[label=\roman*)]
    \setlength\itemsep{0pt} 
        \item if $t^4 r_t^{2d+2}\to 0$, then $\lim_{t\to \infty} \PP(\cross(L_i)=0 \text{ for some } i=1,\ldots,m) =1$.
        \item if $t^4 r_t^{2d+2}\to \infty$, then $\lim_{t \to \infty} \PP(\cross(L_i)=0 \text{ for some } i=1,\ldots,m) = 0$.
        \item if $t^4 r_t^{2d+2}\to c\in(0,\infty)$, then 
    \end{enumerate}
    \[\lim_{t \to \infty} \PP(\cross(L_i)=0 \text{ for some } i=1,\ldots,m) =1-\l(1-M\r)^m\]
    where 
    \[M\coloneq \lim\limits_{t\to \infty} \intl_{\bf G_2(d)} e^{- \EE \cross(L)} \nu_d(\dd L)\] 
    is a constant depending on the limit $c$ and the body $W$. 
\end{theorem}

\begin{remark*} The sampling of planes can be done as follows: 
\begin{itemize}
    \setlength\itemsep{0pt} 
    \item[-] in $\bf G_2(3)$, a uniform plane is given by $u^\perp$ if $u\in S^2$ is uniform;
    \item[-] in arbitrary $\bf G_2(d)$, a uniform plane is given by $\spanvec\{v_1,v_2\}$ if $v_1\in S^{d-1}$ and $v_2\in S^{d-1}\cap v_1^\perp$ are uniform,
\end{itemize} 
where $S^{d-1}$ denotes the boundary sphere of the $d$-dimensional unit ball $B^d$, and $u^\perp$ is the orthogonal complement of $u$.
\end{remark*}

The rest of the paper is structured as follows. In the rest of the introduction, we recall some standard definitions and notation. \Cref{sec:Poisson} contains a Poisson approximation result for the total number of crossings on multiple planes. To prove this result, we also give bounds for the covariance between the number of crossings onto two planes. Section~\ref{sec:proofs} contains the proofs of Theorems~\ref{thm:existence} and \ref{thm:finding}.

\subsection{Preliminaries}

Fix a convex body $W\subset \R^d$ ($d\geq 3$). We note that while it is often assumed that $W$ is unit-volume (in particular in \cite{Chimani2018} and \cite{Döring2025}), this is only for convenience, and the results can be extended to any convex body by rescaling.
Let $\eta_t$ be a homogeneous Poisson point process on $W$ with intensity $t>0$. 
By a standard abuse of notation we associate $\eta_t$ with its support, writing $x\in\eta_t$ whenever $\eta_t(\{x\})=1$.

We frequently use the (multivariate) Mecke formula, which states that if $f:(\R^d)^k\to \R$ is a nonnegative function, then
$$\EE \sum_{(x_1,\ldots, x_k) \in \eta_{t,\neq}^k} f(x_1,\ldots, x_k) = t^k \int\limits_{W^k} \EE f(x_1,\ldots, x_k) \dd x_1 \ldots \dd x_k,$$
where $\eta_{t,\neq}^k$ denotes the set of distinct $k$-tuples of $\eta_t$.

For two real-valued random variables $X$ and $Y$, their total variation distance $\dd_{TV}$ is given by
\[\dTV(X,Y) \coloneq \sup_{A\in\mathcal B(\R)} |\PP(X\in A)-\PP(Y\in A)|.\]
In particular, if $\dTV(X_n,X)\to 0$ as $n\to \infty$ for some seqence $(X_n)_{n\in \N}$ of random variables, convergence in distribution to $X$ follows.
 

We denote the $d$-dimensional Lebesgue measure by $\lambda_d$, and in particular write $\kappa_d\coloneq \lambda_d(B^d)$ for the volume of the ($d$-dimensional) ball with unit radius. Further, recall that $\bf G_2(d)$ denotes the space of two-dimensional linear subspaces in $\R^d$, and $\nu_d$ the associated probability measure. For a plane $L\in \bf G_2(d)$, the orthogonal complement of $L$ is denoted by $L^{\perp}$.

Let $L_1,L_2 \in \bf G_2(d)$. By a classical result of Jordan \cite{Jordan}, an appropriate coordinate system can be chosen such that 
\begin{align}
\begin{split}\label{eq:principalangle}
    L_1 & =\spanvec\{e_1,e_2\} \text{ and}\\
    L_2 &=\spanvec\{e_1 \cos \theta_1+e_3\sin \theta_1, e_2\cos \theta_2 + e_4\sin \theta_2\},
\end{split}
\end{align}
with $e_i$ denoting the $i$th standard unit vector.
The angles $0\leq \theta_1\leq \theta_2 \leq \pi/2$ are referred to as the \textit{principal angles} between the two planes. Note that $\theta_2=\theta_1=0$ is equivalent to $L_1=L_2$, and $\theta_2>\theta_1=0$ to $L_1$ and $L_2$ intersecting in a line, i.e. $L_1\cup L_2$ spanning a three-dimensional linear subspace. In particular, if $d=3$, we always have $\theta_1=0$, while in $d>3$, $\theta_1$ is positive $\nu_d\times \nu_d$-almost everywhere.

Throughout the paper, we write $f(x) = \c O(g(x))$ as $x\to\infty$ if there exist $M>0$ and $x_0$ such that $|f(x)| \leq M|g(x)|$ for all $x\geq x_0$. We also write $f(x) = \c O_x(g(x))$ and omit the $x\to\infty$. Similarly, we write $f(x)=o_x(g(x))$ if for all $M>0$ there exists an $x_0$ such that $|f(x)|\leq M|g(x)|$ for all $x\geq x_0$, or equivalently, if $g(x)\neq 0$, we have $f(x)/g(x)\to 0$ as $x\to \infty$.

\section{Poisson approximation for crossings on multiple planes}\label{sec:Poisson}

To prove the main theorems in this paper, we need to control the dependence between the numbers of crossings $\cross (L_1)$ and $\cross(L_2)$ in projections onto the planes $L_1$ and $L_2$. In this section, we show that if the covariance between $\cross(L_1)$ and $\cross(L_2)$ decreases quickly enough with increasing $t$, then the total number of crossings on these two planes is asymptotically Poisson distributed. 
In particular, if the principal angles between the two planes are large, then the total variation distance to a Poisson random variable is small. This relationship is quantified by the following proposition. 

\begin{proposition}\label{prop:Poisson}
Let $L_1,\ldots, L_k\in \bf G_2(d)$ be arbitrary planes, and denote by $0\leq \theta_1^{i,j}\leq\theta_2^{i,j}\leq \pi/2$ the principle angles between $L_i$ and $L_j$ $(i<j)$. Further, let $\lambda \coloneq \EE \cross L_1$ and $\Lambda = \EE \cross(L_1)+\ldots + \EE \cross(L_k)$, and let $Z_\alpha$ denote a Poisson distributed random variable with expectation $\alpha>0$. Then
\begin{align*}
    \dTV(\cross(L_1),Z_{\lambda}) & \leq D_t, \, \text{ and} \notag \\ 
    \dTV(\cross(L_1) + \ldots + \cross(L_k), Z_{\Lambda}) & \leq kD_t +  c t^4 r_t^{2d} \sum_{1\leq i<j\leq k} b_{\theta_1^{i,j}, \theta_2^{i,j}}, 
    \end{align*}
     where $c>0$ is a constant,
        \[D_t = 
        \begin{cases}
            \c O_t(t^7 r_t^{4d+4}),& \text{if }  \lim_{t\to \infty} tr_t^d >0; \\
            \c O_t(t^6r_t^{3d+4}), & \text{if } \lim_{t\to \infty} tr_t^d = 0 \text{ and } \lim_{t\to \infty} t r_t^2>0;\\
            \c O_t(t^4r_t^{3d}), & \text{if }\lim_{t\to \infty} t r_t^2 = 0,
        \end{cases}\]
        and 
        \[b_{\theta_1,\theta_2} = \min\l\{r_t^2,\ r_t^3 (\sin \theta_2)^{-1},\ r_t^4 (\sin \theta_1 \sin \theta_2)^{-1}\r\},\]
        where the implied constants in $D_t$ may depend on $d$ and $W$, but not on $L_1$ or $L_2$.
\end{proposition}

The proof of the proposition relies on the following lemma.

\begin{lemma}\label{lemma:covbound}
    Let $L_1, L_2\in \bf G_2(d)$, and let $0\leq \theta_1\leq\theta_2\leq \pi/2$ denote the principle angles between them. Then
    \begin{align*}
        \VV(\cross(L_1)) & = D_t + \EE (\cross(L_1)), \\  
        \Cov(\cross(L_1), \cross(L_2)) &= D_t + \c O_t(t^4 r_t^{2d} \cdot b_{\theta_1,\theta_2}),
    \end{align*}
    with $D_t$ and $b_{\theta_1,\theta_2}$ given in \Cref{prop:Poisson}.
\end{lemma}

We first show how this result is used to prove the proposition above.

\begin{proof}[Proof of \Cref{prop:Poisson}]
    We use the Poisson approximation result \cite[Theorem~3.8]{Pianoforte} (see in particular (3.27)), which states the following: if $S$ is a Poisson $U$-statistic of order $l$ with symmetric kernel function, and $Z\sim \Poi(\EE S)$, then the total variation distance is upper bounded by
    \begin{equation*}
        \dTV(S,Z)\leq 2^l \cdot \min\l(1, \EE S\r) \cdot \frac{\VV S - \EE S}{\EE S}\leq 2^l \cdot (\VV S-\EE S)
    \end{equation*}
    Since $\cross (L_1)$, as well as the sum, is a Poisson $U$-statistic (of order $4$), this result can be applied. For a single plane, $\dTV(\cross(L_1),Z_{\lambda}) = \c O_t(\VV \cross(L_1)-\EE \cross(L_1))$, and for the sum,
    \[\dTV(\cross(L_1) + \ldots + \cross(L_k), Z_{\Lambda}) = \sum_{i=1}^k\c O_t \big(\VV \cross(L_i)-\EE \cross(L_i)\big) + \sum_{1\leq i<j\leq k} \c O_t\big( \Cov(\cross(L_i),\cross(L_j))\big).\]
    The statement of the proposition then follows from \Cref{lemma:covbound}.
\end{proof}

    Note that the statement for the variance has already been shown in \cite[Section~4]{Chimani2018}, see also \cite[Lemma~2.2]{Döring2025}. 

    \begin{proof}[Proof of \Cref{lemma:covbound}.]
    We aim to give an upper bound on $\EE[\cross(L_1) \cross(L_2)]-\EE \cross(L_1) \EE \cross(L_2)$. Unless stated otherwise, we allow $L_1=L_2$, leading to the term $D_t$ that the variance and covariance share.
    
    For the expectation of the product, we have
    \[\EE [\cross(L_1)\cross(L_2)] = \frac{1}{24^2}\cdot \EE \sum_{\substack{(v_1,\dots,v_4)\in \eta_{t,\neq}^4 \\ (w_1,\dots,w_4)\in \eta_{t,\neq}^4}} h_{L_1}(v_1,v_2,v_3,v_4) h_{L_2}(w_1,w_2,w_3,w_4),\]
    where $h_L$ denotes the indicator function that the points generate a crossing when projected onto $L$, defined precisely in (\ref{eq:kernel}). 
    If all vertices in the sum are distinct, then an application of the Mecke formula results in the product of the expectations, $\EE \cross(L_1) \EE \cross(L_2)$. 
    As a consequence, the covariance is given by the evaluation of the above sum over $\{v_1,\ldots, v_4\}$ and $\{w_1,\ldots, w_4\}$ with nonempty intersection. 
    Note that as we are only interested in the order of magnitude, we omit constants, and simply write a generic constant $c$, which might change from line to line; however, it is always independent of $t$ and universal in the planes, but depends on $d$ and $W$. Further, for simplicity, we only consider the case when $W=R\cdot B^d$; any such upper bound is also an upper bound for arbitrary convex body $W'$ with circumradius at most $R$.

    To apply the Mecke formula again, we have to distinguish cases based on the cardinality $k$ of the intersection $\{v_1,v_2,v_3,v_4\}\cap \{w_1,w_2,w_3,w_4\}$; denote these sums by $I_k$ for $k=1,\ldots, 4$. 

    For $k=1$, we have
    \begin{align*}
        I_1\coloneq & c \EE \sum_{(v_1,v_2,v_3,v_4,w_2,w_3,w_4) \in \eta_{t,\neq}^7} h_{L_1}(v_1,v_2,v_3,v_4) h_{L_2}(v_1,w_2,w_3,w_4) \\
        = & c t^7\int_{W^7} 
           \begin{multlined}[t]
                 \ind(\|v_1-v_2\|\leq \radt,\, \|v_3-v_4\|\leq \radt,\, \|v_1-w_2\|\leq \radt,\, \|w_3-w_4\|\leq \radt)\\
                \times \ind([v_1,v_2]|_{L_1}\cap [v_3,v_4]|_{L_1}\neq \emptyset, \, [v_1,w_2]|_{L_2}\cap [w_3,w_4]|_{L_2}\neq \emptyset) \dd v_1 \cdots \dd v_4 \dd w_2 \dd w_3 \dd w_4.
            \end{multlined}.
    \end{align*}
    Note that as both $h_{L_1}$ and $h_{L_2}$ are written as a sum, the integrand should be a sum as well; however, it is easy to see that every summand actually yields the same integral, only contributing a constant factor to the expression.

    Now, the crossing of two edges implies that the projection of any two of the four points have distance at most $2r_t$ (on the corresponding plane), hence    
    \begin{multline*}
        I_1 
        \leq ct^7\int_{W^7} \ind(\|v_1-v_2\|\leq \radt,\, \|v_3-v_4\|\leq \radt,\, \|v_1-w_2\|\leq \radt,\, \|w_3-w_4\|\leq \radt) \\
        \times \ind(\|v_3|_{L_1} - v_1|_{L_1} \|\leq 2\radt,\, \|w_3|_{L_2} - v_1|_{L_2}\|\leq 2\radt) \dd v_1 \cdots \dd v_4 \dd w_2 \dd w_3 \dd w_4,
    \end{multline*}
    For fixed $v_1, v_3, w_3$, the remaining $4$ vertices are each contained in a ball of radius $r_t$ (around the other endpoints of their respective edges), leading to the upper bound 
    \[I_1\leq ct^7 (r_t^d)^4\int_{W^3} \ind\l(\|v_3|_{L_1} - v_1|_{L_1} \|\leq 2\radt,\, \|w_3|_{L_2} - v_1|_{L_2}\|\leq 2\radt\r) \dd v_1 \dd v_3 \dd w_3.\]
    Further, we have that
    \begin{equation*}
        \sup_{v_1\in W} \lambda_d\l(\l\{ v_3\in W\colon\, \|v_3|_{L_1}-v_1|_{L_1}\|\leq 2\radt\r\}\r) \leq 2R \cdot (2 \radt)^2\kappa_2,
    \end{equation*}
    where $R$ denotes the radius of the ball $W$.
    Applying this to $v_3$ and $w_3$, and using that $v_1$ is contained in $W$, this leads to an upper bound on $I_1$ of order $t^7 r_t^{4d} (r_t^2)^2 = t^7 r_t^{4d+4}$.

    For $k=2$, we consider
    \[I_2 = ct^6\int_{W^6} 
       h_{L_1}(v_1,v_2,v_3,v_4) \cdot h_{L_2}(v_1,v_2,w_3,w_4)\dd v_1 \cdots \dd v_4\dd w_3 \dd w_4.\]

    In this case, we have to differentiate three geometric cases based on which pair of edges produces the crossing in each tuple: either the coinciding pair $(v_1,v_2)$ is an edge involved in the forming of the crossing on both planes, only on one plane, or neither. Let us denote these integrals by $I_2^{(1)}$, $I_2^{(2)}$, and $I_2^{(3)}$, respectively.
    
    The argument of the bounds follow similarly as above: keeping a single vertex fixed, the rest of the vertices are either constrained by the crossings (giving a distance constraint of $2r_t$ on the appropriate plane, which generates a factor of $r_t^2$), or by the existence of an edge (generating a factor of $r_t^d$).
    It thus follows that
    \begin{align*}
        I_2^{(1)} 
        = & c t^6\int_{W^6} 
           \begin{multlined}[t]
                 \ind(\|v_1-v_2\|\leq \radt,\, \|v_3-v_4\|\leq \radt,\, \|w_3-w_4\|\leq \radt)\\
                 \times \ind([v_1,v_2]|_{L_1}\cap [v_3,v_4]|_{L_1}\neq \emptyset, \, [v_1,v_2]|_{L_2}\cap [w_3,w_4]|_{L_2}\neq \emptyset) \dd v_1 \cdots \dd v_4 \dd w_3 \dd w_4
            \end{multlined}\\ 
        & \leq ct^6 \int_{W^6}  \begin{multlined}[t]
        \ind \l(\|v_1-v_2\|\leq r_t, \|v_3-v_4\|\leq r_t, \|w_3-w_4\|\leq r_t)\r) \\ \times \ind\l( \|v_1|_{L_1} - v_3|_{L_1}\|\leq 2r_t, \|v_1|_{L_2} - w_3|_{L_2}\|\leq 2r_t\r) \dd v_1 \cdots \dd v_4 \dd w_3 \dd w_4
            \end{multlined}\\
        & \leq c t^6 (r_t^d)^3 (r_t^2)^2 = c t^6 r_t^{3d+4}
    \end{align*}
    with $v_3$ and $w_3$ being constrained by the existence of a crossing, and $v_2$, $v_4$ and $w_4$ by the existence of an edge. 
    Similarly, we  have
    \begin{align*}
        I_2^{(2)}
        & \leq c t^6 \int_{W^6}  \begin{multlined}[t]
        \ind \l(\|v_1-v_2\|\leq r_t, \|v_3-v_4\|\leq r_t, \|v_1-w_3\|\leq r_t, \|v_2-w_4\|\leq r_t\r) \\ \ind\l( \|v_1|_{L_1} - v_3|_{L_1}\|\leq 2r_t\r) \dd v_1 \cdots \dd v_4 \dd w_3 \dd w_4 \leq
            \end{multlined}\\
        & \leq c t^6 (r_t^d)^4 (r_t^2) = c t^6 r_t^{4d+2}.
    \end{align*}
    and
    \begin{align*}
        I_2^{(3)}
        & \leq ct^6 \int_{W^6}  \begin{multlined}[t]
        \ind \l(\|v_1-v_3\|\leq r_t, \|v_2-v_4\|\leq r_t, \|v_1-w_2\|\leq r_t, \|v_3-w_4\|\r) \\ \ind\l( \|v_1|_{L_1} - v_2|_{L_1}\|\leq 2r_t\r) \dd v_1 \cdots \dd v_4 \dd w_3 \dd w_4 \leq
            \end{multlined}\\
        & \leq c t^6 (r_t^d)^4 (r_t^2) = ct^6 r_t^{4d+2}.
    \end{align*}
    Note that in the two latter cases the second crossing provides no additional condition on the upper bound; this is because the closeness of the vertices in the projection is guaranteed by the existence of the edges themselves.

    Now, let $k=3$, say $v_i=w_i$ for $i=1,2,3$. For both planes, there is a pair in the coinciding triple that is an edge in the crossing; this pair either coincides (i.e. the same edge produces the crossing on both planes), or it does not.
    Denoting these by $I_3^{(1)}$ and $I_3^{(2)}$, respectively, we have 
     \begin{align*}
         I_3^{(1)} 
         & \leq ct^5 \int_{W^5}  \begin{multlined}[t]
        \ind \l(\|v_1-v_2\|\leq r_t, \|v_3-v_4\|\leq r_t, \|v_3-w_4\|\leq r_t\r) \\ \times \ind\l( \|v_1|_{L_1} - v_3|_{L_1}\|\leq r_t\r) \dd v_1 \cdots \dd v_4 \dd w_3 \dd w_4 \leq
             \end{multlined}\\
         & \leq c t^5 (r_t^d)^3 r_t^2 = c t^5 r_t^{3d+2}.
     \end{align*}
    and
    \begin{align*}
        I_3^{(2)} 
        & \leq ct^5 \intl_{W^5}
        \begin{multlined}[t]
        \ind \l(\|v_1-v_2\|\leq r_t,  \|v_3-v_4\|\leq r_t\r) \\ \times \ind\l( \|v_1-v_3\|\leq r_t, \|v_2-w_4\|\leq r_t\r)\dd v_1 \cdots \dd v_4 \dd w_4 \leq 
             \end{multlined}\\
        & \leq c t^5 (r_t^d)^4= c t^5 r_t^{4d}.
    \end{align*}


    Lastly, we consider the case $k=4$. We again have two possibilities: either the crossing-forming edges coincide, or they do not. In the latter case, any two points have distance at most $2r_t$, and thus the corresponding integral has order of magnitude $t^4 r_t^{3d}$. In the former case, we have that there exist two edges, and they form a crossing on both planes
    Thus, the integral, denoted by $I_4$, becomes
    \begin{multline*}
        I_4=ct^4\int_{W^4} \ind \l(\|v_1-v_2\|\leq \radt,\, \|v_3-v_4\|\leq \radt \r) \\
        \times \ind \l([v_1,v_2]|_{L_1}\cap [v_3,v_4]|_{L_1}\neq \emptyset , [v_1,v_2]|_{L_2}\cap [v_3,v_4]|_{L_2}\neq \emptyset \r) \dd v_1 \dd v_2 \dd v_3 \dd v_4.
    \end{multline*}
    If $L_1=L_2$, this expression is equal to the expectation, yielding the claim of the lemma for the variance. Otherwise, we upper bound the expression as before by
    \[I_4\leq ct^4 r_t^{2d} \int_{W^2} \ind\l( \|v_1|_{L_1} - v_3|_{L_1}\|\leq 2r_t, \|v_1|_{L_2} - v_3|_{L_2}\|\leq 2r_t\r) \dd v_1 \dd v_3.\]
    The condition $\|v|_L-w|_L\|\leq \delta$ is equivalent to $w|_L\in v|_L+\delta B^d|_L$, and thus $w\in v|_L+\delta B^d + L^\perp = v+\delta B^d +L^\perp$. Writing $P_L^{\delta}\coloneqq L^\perp + \delta B^d$ then gives
    \[I_4\leq ct^4 r_t^{2d} \lambda_d(P_{L_1}^{2r_t} \cap P_{L_2}^{2r_t}\cap (2W)).\]

    To conclude the proof, we show that $\lambda_d(P_{L_1}^{2r_t} \cap P_{L_2}^{2r_t}\cap (2W)) \leq b_{\theta_1,\theta_2}$ holds.
    This volume depends on the closeness of the planes, or in particular, on their principal angles. 
    We use the coordinate system in which (\ref{eq:principalangle}) holds, which is allowed by the symmetry of $W=RB^d$.
    Denoting the spanning vectors of $L_2$ by 
    \[u_1\coloneq e_1\cos \theta_1+e_3 \sin \theta_1 \text{ and } u_2\coloneq e_2\cos \theta_2 + e_4 \sin \theta_2,\]
    we obtain
    \[P_{L_1}^{2r_t} \cap P_{L_2}^{2r_t}\cap (2W) \subset \{x\in 2W : |\langle x,e_i\rangle|,|\langle x,u_i\rangle|\leq 2r_t \text{ for } i=1,2\}.\]

    The last $d-4$ coordinates are only constrained through the containment in $2W$, independently of $r_t$. For the first $4$ coordinates, we obtain an upper bound by dropping the condition of being contained in $2W$, and considering only the inner products. By a standard determinant computation, this gives
    \[\lambda_d(P_{L_1}^{2r_t} \cap P_{L_2}^{2r_t}\cap (2W))\leq c r_t^4 (\sin \theta_1 \sin \theta_2)^{-1}.\]

    If $\theta_1$ is small (and in particular, if $\theta_1=0$, which is always the case if $d=3$), the containment in $2W$ is stronger than the inner product condition for the third coordinate. Dropping the inner product condition with $u_1$ results in the upper bound
     \[\lambda_d(P_{L_1}^{2r_t} \cap P_{L_2}^{2r_t}\cap (2W))\leq c r_t^3 (\sin \theta_2)^{-1}.\]
     Lastly, if both angles are very small, the inner products with $u_i$-s only give weak constraints. In particular, if we drop the condition with $u_1$ and $u_2$, we obtain the bound $\lambda_d(P_{L_1}^{2r_t} \cap P_{L_2}^{2r_t}\cap (2W))\leq cr_t^2$.
    \end{proof}

\section{Proof of theorems}\label{sec:proofs}

\subsection*{Existence of a plane without crossings}\label{sec:existence}
\begin{proof}[Proof of \Cref{thm:existence}]

We only prove the statement for $W=R\cdot B^d$. Restricting the model to an arbitrary subset of the ball can only remove edges, and thus potential crossings, increasing the probability in question. The claim for general $W$ easily follows. 
In addition, fix an arbitrary plane $L_0$. By symmetry, every statement for $L_0$ holds for any arbitrary plane.

By the second moment method, we can bound the probability that there
exists a plane without crossings from below:
    \begin{equation*}
        \PP\l(X_t >0\r) \geq \frac{(\EE X_t)^2}{\EE X_t^2}.
    \end{equation*}
    We have
    \[ \EE X_t = \int_{\bf G_2(d)} \PP(\cross(L)=0) \nu_d(\dd L) = \PP(\cross(L_0)=0)\]
    and
    \begin{align*}
        \EE X_t^2 
        & = \int_{\bf G_2(d)}\int_{\bf G_2(d)} \PP\l(\cross(L_1) = 0,\,\cross(L_2)  = 0\r) \nu_d(\dd L_1) \nu_d(\dd L_2) \\
        & = \int_{\bf G_2(d)}\int_{\bf G_2(d)} \PP\l(\cross(L_1) + \cross(L_2) = 0\r) \nu_d(\dd L_1) \nu_d(\dd L_2) \\
        & = \int_{\bf G_2(d)} \PP\l(\cross(L_0) + \cross(L_1) = 0\r) \nu_d(\dd L_1).
    \end{align*}

    To bound the probabilities appearing in the integral, we use \Cref{prop:Poisson}.
    We assume that $tr_t^d\to 0$ (which gives an upper bound on $r_t$ that shall be fulfilled by our final condition), and for simplicity, that $t^2r_t^{d+1}\to \infty$ (or equivalently, $\lambda \to \infty$). The latter gives a lower bound on $r_t$; however, if the statement of the theorem holds for a certain threshold distance, by monotonicity it clearly holds for smaller ones, and so this is essentially not a restriction.
    Since $d\geq 3$, the latter condition implies that $t r_t^2$ is bounded away from zero, and $D_t = \c O_t(t^6 r_t^{3d+4})$ in \Cref{prop:Poisson}.

    From \cite[Lemma~1]{Döring2025}, we have that 
    \[\lambda \coloneq \EE \cross(L_0) = c_Rt^4 r_t^{2d+2} + \c O_t(r_t),\]
    with $c_R$ being dependent on $R$ and $d$.
    For the first moment, \Cref{prop:Poisson} immediately gives
    \[\EE X_t = \PP(\cross(L_0)=0) = e^{-\lambda} + \c O_t(t^6 r_t^{3d+4}).\]

    For the second moment, we handle the pairs of planes based on their covariance.
    Set
    \begin{align}\label{eq:closeplanes}
        C_\delta \coloneq  
        \begin{cases}
            \l\{L_1 \in \bf G_2(d):\sin \theta_2<\delta \r\}, & \text{ if } d=3, \\
            \l\{L_1 \in \bf G_2(d):\sin \theta_1 \sin \theta_2<\delta \r\},&  \text{ if } d\geq 4,
        \end{cases}
    \end{align}
    where $\theta_1,\theta_2$ are the principal angles between $L_0$ and $L_1$.
    This set corresponds to planes having a small angle, and thus large covariance, with $L_0$.
    For arbitrary $\delta>0$, we obtain the upper bound
    \begin{align*}
        \EE X_t^2 
        & \leq \intl_{C_\delta}\PP\l(\cross(L_0)= 0\r) \nu_d(\dd L_1) +\intl_{\bf G_2(d) \setminus C_\delta} \PP\l(\cross(L_0) + \cross(L_1) = 0\r) \nu_d(\dd L_1)  \\
        & \leq \nu_d(C_\delta)\PP\l(\cross(L_0) = 0\r) + \sup_{L\in \bf G_2(d)\setminus C_\delta} \PP\l(\cross(L_0) + \cross(L) = 0\r).
    \end{align*}

    For the second term, we have from \Cref{prop:Poisson} that for all $L\in \bf G_2(d)\setminus C_{\delta}$,
    \[\PP\l(\cross(L_0) + \cross(L) = 0\r) = e^{-2\lambda} + \c O_t(t^6 r_t^{3d+4}) +\c O_t \l(t^4 r_t^{2d+\min\{4,d\}} \delta^{-1}\r).\]
    For the former, we need to find the measure of planes having large covariance with $L_0$.

    \begin{claim}\label{claim:closeplanes}
        With $C_\delta$ defined by (\ref{eq:closeplanes}), it holds that
        \[\nu_d(C_\delta)=\begin{cases}
            \c O_t(\delta^{d-3} \ln(\delta^{-1})), & \text{ if } d\geq 4,\\
            \c O_t(\delta^{2}), & \text{ if } d=3 
        \end{cases}\]
        as $\delta \to 0$.
    \end{claim}
    \begin{proof}
    An integral with respect to $\nu_d$ can be transformed to one with respect to the principal angles: see \cite[Section~7]{James}.
    For $d=3$, it follows from \cite[(7.17)]{James} that
    \[\nu_3(C_\delta) = \int_{\bf G_2(3)} \ind ( \sin \theta_2<\delta)\nu_3(\dd L) =\int_0^{\pi/2} \ind\l( \sin \theta_2 <\delta\r)\sin \theta_2 \dd \theta_2.\]
    Explicit integral computation gives that $\nu_3(C_\delta)$ has order of magnitude $\delta^2$ as $\delta\to 0$. 

    For $d\geq 4$, we obtain 
    \begin{align*}
        \nu_d(C_\delta)
        & =\int_{\bf G_2(d)}\ind\l( \sin \theta_1 \sin \theta_2 <\delta\r) \nu_d(\dd L) \\
        & = c \iint\limits_{0\leq \theta_1 \leq \theta_2 \leq \frac{\pi}{2}} \ind(\sin \theta_1 \sin \theta_2<\delta) (\sin \theta_1 \sin \theta_2)^{d-4} (\sin^2(\theta_2)- \sin^2(\theta_1))\dd \theta_1 \dd \theta_2
    \end{align*}
    from \cite[(7.13)]{James}.
    We use the bound given by the indicator for the product of sines, and the trivial upper bound of $1$ for the difference of squares. We are then left with only an indicator, and the inner integral can be computed explicitly: 
    \begin{align*}
        \nu_d(C_\delta) 
        & \leq c\delta^{d-4} \int_{0}^{\frac{\pi}{2}} \int_0^{\frac{\pi}{2}} \ind(\sin \theta_1 \sin \theta_2 < \delta)\dd \theta_1 \dd \theta_2 \\
        & \leq c\delta^{d-4}\l[\int_0^{\arcsin \delta} \dd \theta_2 + \int_{\arcsin \delta}^{\frac{\pi}{2}} \arcsin\l(\frac{\delta}{\sin \theta_2}\r)\dd \theta_2\r].
    \end{align*}
    We use the crude analytical upper bounds $\arcsin \l(\frac{\delta}{\sin \theta_2}\r)\leq \frac{2\delta}{\sin \theta_2} \leq \frac{4\delta}{\theta_2}$ for the second integral, and obtain
    \[\nu_d(C_\delta)
        \leq c\delta^{d-3}\l[\arcsin \delta + \ln \frac{\pi}{2} - \ln(\arcsin \delta) \r],\]
    which is dominated by the last term, as $-\ln(\arcsin \delta) = \ln\l(\frac{1}{\arcsin \delta}\r)= c\ln\l(\frac{1}{\delta}\r)$, which was the upper bound stated in the claim.
    \end{proof}

    Plugging the asymptotics of the measure of close planes into the second moment method lower bound, we obtain for $d\geq 4$ that
    \[\frac{(\EE X_t)^2}{\EE X_t^2} = \frac{e^{-2\lambda}+ \c O_t(e^{-\lambda} \cdot t^6r_t^{3d+4}) + \c O_t\l(t^{12} r_t^{2(3d+4)}\r)}{e^{-2\lambda} + \c O_t(t^6r_t^{3d+4}) + \c O_t\l(t^4r_t^{2d+4} \cdot \delta^{-1}\r) + \c O_t \l(\delta^{d-3} \ln(\delta^{-1})e^{-\lambda}\r)},\]
    where $\delta=\delta_t\to 0$ as $t\to \infty$ is assumed in the $\c O_t(\cdot)$ notation. We aim to choose $\lambda$ and $\delta$ such that the first exponential expression dominates both numerator and denominator; recall that $\lambda$ is of the order $t^4 r_t^{2d+2}$.
    If $\lambda \leq c\ln t$, then $e^{-2\lambda}\geq t^{-2c}$, and
    \[t^6 r_t^{3d+4} = \c O_t\l(\l(\ln t\r)^{\frac{3d+4}{2d+2}} t^{-\frac{2}{d+1}}\r) \text{ and } t^4 r_t^{2d+4} = \c O_t\l(\l(\ln t\r)^{\frac{2d+4}{2d+2}} t^{-\frac{4}{d+1}}\r).\]
    Comparing the first term to $e^{-2\lambda}$, we immediately see that $c<2/(d+1)$ is necessary and sufficient for $e^{-2\lambda}$ to dominate the numerator.
    Further, if we set $\delta=\delta_t=e^{-c'\lambda}$, with 
    \[\frac{1}{d-3} < c' < 2(d-3),\]
    then $e^{-2\lambda}$ dominates both $\delta^{d-3} \ln (\delta^{-1}) e^{-\lambda}$ and $t^4 r_t^{2d+4} \delta^{-1}$, and thus the entire denominator. Further, as we assumed $\lambda\to \infty$, we also have $\delta\to 0$. As a consequence, the quotient tends to $1$ as $t\to \infty$ whenever $\lambda \leq c\ln t$ with $c<2/(d+1)$.
    For $d=3$, only the last term of the denominator changes, giving different bounds for $c'$, but the same requirement for $c$, which concludes the proof of the theorem.
    \end{proof}

\subsection*{Finding a plane without crossings}\label{sec:finding}

\begin{proof}[Proof of \Cref{thm:finding}.]
Recall that $\EE \cross(L)$ is of the order $t^4 r_t^{2d+2}$ for any fixed plane $L$.

If $t^4 r_t^{2d+2}\to 0$, then 
\[\PP(\cross(L_1)=0) = 1- \PP(\cross(L_1)\geq 1)\geq 1- \EE \cross(L_1)\]
by Markov's inequality, and since the expectation is exactly of the order above, the lower bound converges to $1$. If $t^4 r_t^{2d+2}\to \infty$,
\[\PP(\cross(L_i)=0 \text{ for some } i=1,\ldots, m)\leq m \PP(\cross(L_1)=0)\leq \frac{m\VV \cross(L_1)}{(\EE \cross(L_1))^2},\]
which then tends to $0$ by the asymptotic behaviour of the variance (see \Cref{lemma:covbound}).

Assume now that $t^4 r_t^{2d+2}\to c\in(0,\infty)$. 
We find the limit of $\PP(\cross(L_1)+\ldots + \cross(L_k)=0)$ for arbitrary $k$, from which the statement of the theorem follows via the inclusion-exclusion principle.

For fixed (deterministic) planes $L_1',\ldots, L_k'\in \bf G_2(d)$ with pairwise principal angles $0\leq \theta_1^{i,j}\leq\theta_2^{i,j}\leq \pi/2$ , we have by \Cref{prop:Poisson} that 
\[\dTV(\cross(L_1') + \ldots + \cross(L_k'), Z_{\Lambda}) \leq \c O_t(t^6r_t^{3d+4}) + ct^4 r_t^{2d} \sum_{1\leq i<j\leq k} b_{\theta_1^{i,j}, \theta_2^{i,j}},\]
with $b_{\theta_1,\theta_2}$ given in \Cref{lemma:covbound}, and $\Lambda \coloneq \EE (\cross(L_1') + \ldots + \cross(L_k'))$.
Let $\delta = \delta_t \coloneq t^{-\frac{1}{d+1}}$, and 
\[\c C_\delta = \begin{cases}
    \{(L_1',\ldots, L_k')\in (\bf G_2(d))^k\colon \sin \theta_2^{i,j} <\delta\}, & \text{ if } d=3, \\
    \{(L_1',\ldots, L_k')\in (\bf G_2(d))^k \colon \sin \theta_1^{i,j} \sin \theta_2^{i,j} <\delta\}, & \text{ if } d\geq 4.
\end{cases}\]
Since the planes are chosen uniformly, we have that
\[\PP(\cross (L_1)+\ldots + \cross(L_k) = 0 )  = \idotsint\limits_{(\bf G_2(d))^k} \PP(\cross(L_1')+\ldots + \cross (L_k')=0) \nu_d(\dd L_1')\cdots \nu_d(\dd L_k'),\]
which we divide into integrals over $\c C_\delta$ and its complement. 
Then on one hand, the total variation distance above tends to $0$ for all planes in $\c C_\delta$, and convergence in distribution to the Poisson random variable $Z_{\Lambda}$ follows. On the other hand, $\nu_d^k(\c C_\delta)$ tends to $0$ by \Cref{claim:closeplanes}, which yields
\begin{multline*}
 \lim_{t\to \infty} \PP(\cross(L_1)+\ldots +\cross(L_k)=0) \\ 
 = \idotsint\limits_{(\bf G_2(d))^k} e^{-(\lambda_{L_1'}+\ldots + \lambda_{L_k'})} \nu_d(\dd L_1')\cdots \nu_d(\dd L_k') 
 = \l[ \intl_{\bf G_2(d)} e^{- \lambda_{L}} \nu_d(\dd L)\r]^k,
\end{multline*}
where $\lambda_L\coloneq \lim_{t\to \infty} \EE \cross(L)$. 
Writing $M$ for the last integral, we have from the inclusion-exclusion principle that
\begin{align*}
    \PP(\cross(L_i)=0 \text{ for some } i=1,\ldots, m) 
    & = \sum_{n=1}^m (-1)^{n-1} \binom mn \PP(\cross(L_1)=\ldots = \cross(L_n)=0) \\
    & = \sum_{n=1}^m (-1)^{n-1}\binom mn \PP(\cross(L_1)+\ldots + \cross (L_n)=0)
\end{align*}
hence by the previous limit and the binomial expansion, 
\[\lim_{t\to \infty} \PP(\cross(L_i)=0 \text{ for some } i=1,\ldots, m) = \sum_{n=1}^m (-1)^{n-1}\binom mn M^n = 1-(1-M)^m \]
concludes the proof of the theorem.
\end{proof}

\paragraph{Acknowledgements}
The authors would like to thank Júlia Komjáthy for suggesting this problem at the Workshop On Randomness and Discrete Structures in Groningen.
We also thank Bernhard Hafer for pointing out a mistake in a previous version of the manuscript.

\paragraph{Funding information}
This research was funded by the Deutsche Forschungsgemeinschaft (DFG, German Research Foundation) – Project-ID 531542011 (L.~de Jonge) and 531562368 (K.~Nagy). The research was carried out while L.~de~Jonge was employed at Osnabrück University.

\addcontentsline{toc}{section}{References}
\bibliographystyle{abbrv}

\end{document}